\documentclass[ 10pt] {amsart}	
\usepackage{amscd}
\usepackage[all]{xy}

\newtheorem{theo}{Theorem}[section]
\newtheorem{prop}[theo]{Proposition}

\newtheorem{la}[theo]{Lemma }

\theoremstyle{remark}
\newtheorem*{rem}{Remark}
\newtheorem*{note}{Note}

\theoremstyle{definition}

\newtheorem{exmp}{Example}[section]

\newcommand{\N}{\mathbb N}
\newcommand{\C}{\mathbb C}
\newcommand{\R}{\mathbb R}
\newcommand{\Q}{\mathbb Q}
\newcommand{\Z}{\mathbb Z}
\newcommand{\co}{\mathcal O}
\newcommand{\A}{\mathbf A}

\DeclareMathOperator{\reg}{\rm reg}
\DeclareMathOperator{\coker}{\rm coker}
\DeclareMathOperator{\im}{\rm im}
\DeclareMathOperator{\Hom}{\rm Hom}
\DeclareMathOperator{\Ext}{\rm Ext}
\DeclareMathOperator{\Fil}{\rm Fil}
\DeclareMathOperator{\tor}{\rm tor}
 
\numberwithin{equation}{section}
\begin{document}
\title[universal distribution and index formula]{Spectral sequence of 
universal distribution
and Sinnott's index formula}
\author{Yi Ouyang}
\address{
School of Mathematics, University of Minnesota, Minneapolis, 
MN 55455, USA} 
\email{youyang@math.umn.edu}
\date{May 12, 1999}

\begin{abstract}
We prove an abstract index formula about Sinnott's symbol between two
different lattices.   We also develop the theory of the universal distribution and 
predistribution in a double complex point of view. The theory of spectral sequence is used
to interpret the index formula and to analyze the cohomology of  the universal distribution.
Combing these results,   we successfully prove  Sinnott's index formula about the Stickelberger ideal.
In addition, the $\{\pm 1\}$-cohomology groups of  the universal distribution and the 
universal predistribution are obtained.
\end{abstract}
\maketitle
\section{Introduction}
The theory of universal distribution, with its tremendous application in number theory, has
been well studied in the past thirty years(See Lang~\cite{Lang1} and Washington~\cite
{Washington} for more backgrounds).  In \cite{Yamamoto} , Yamamoto studied
the $\{\pm 1\}$-cohomology of the  universal distribution of rank 1(the ${gap}$ group in 
\cite{Yamamoto}):
\begin{theo} \label{ta}
Let $U_m$ be the universal distribution of rank $1 $ and level $m$. Then
\[ H^i(\{\pm 1\}, U_m)=(\Z/2\Z)^{2^{r-1}} \]
where $r$ is the number of distinct prime factors of $m$. \end{theo}
 
In his famous paper
~\cite{S1}, Sinnott successfully obtained the index formula of Stickelberger ideal and circular 
units, which generalized the results of Kummer and Iwasawa. His result can be stated as
\begin{theo} \label{tb}
Let $m$ be a positive integer which is not $2\pmod 4$. Let $G$ be the Galois
group of the cyclotomic extension $\Q(\zeta_m)/\Q$.  Let $R=\Z[G]$ and let $S$ be the Stickelberger
ideal of $\Q(\zeta_m)$. Let $E$ be the group of units in $\Q(\zeta_m)$ and let $C$ be the subgroup
of circular units in $E$. then
 
(1). $ [R^{-}:S^{-}=2^a h^{-}$;

(2). $[E^{+}: C^{+}]=2^b h^{+}$;

\noindent where $a=b=0$ if $r=1$ and $a=2^{r-2}-1$, $b=2^{r-2}+1-r$ if $r>1$,  $h^{+}$
and $h^{-}$ are the class number of $\Q(\zeta_m)^{+}$ and the relative class number of $\Q(\zeta_m)$
respectively.
\end{theo}

Sinnott's result was a huge success and inspired many followers. Most notably, 
Kubert~\cite{Kubert1} and \cite{Kubert2}  found the connection of Sinnott's method
and the universal (ordinary) distribution. By using this connection, he thus
showed that Theorem~\ref{ta} is true for the universal distribution of arbitrary rank.

Sinnott's computation is very elegant but rather  difficult. The motivation to find an easier
proof drives us to this paper. The theory of spectral sequences, though very popular in 
topology, algebra and even number theory, had not been able to leave its mark in the 
theory of distribution until recently. In \cite{ad2}, Anderson came up with 
the  idea  of  using a special double complex to compute the $\{\pm 1\}$-cohomology 
of the universal distribution.  With which he proved a conjecture by Yin~\cite{Yin}. 
Das~\cite{Das}  then used it  to study algebraic monomials and obtained many interesting results.
Their method is  the prototype of spectral sequences method used by us here. 

In this paper,  we prove an abstract index formula about Sinnott's symbol between two
different lattices.   We also develop the theory of the universal distribution and 
predistribution in a double complex point of view. The theory of spectral sequence is used
to interpret the index formula and to analyze the cohomology of  the universal distribution.
Combing these results,   we successfully prove  Sinnott's index formula about the Stickelberger ideal
(i.e., the first part of Theorem~\ref{tb}).
In addition, the $\{\pm 1\}$-cohomology groups of  the universal distribution(i.e.,  Theorem~\ref{ta})
 and the universal predistribution are obtained.  Although we only study the rank $1$ case in this 
paper,  our method is capable of generalizing to the higher rank case. 

As noted above, this paper is based on my advisor,  Professor Greg W. Anderson's
brilliant idea. I am in debt to his working note~\cite{ad1} which contains the raw form
of the abstract index formula and many other facts stated in this paper.
I also benefit greatly from numerous discussions with him.  This paper would be impossible 
without his instruction.  I thank whole heartedly for his insight, patience and encouragement. 

\section{The abstract index formula}
\subsection{Definition of regulator  $\reg(A, B,\lambda)$}
Let $A$ and $B$ be lattices in a finite dimensional 
vector space $V$ over $\R$. Necessarily there exists some $\R$-linear 
automorphism $\phi$ of $V$ such that $\phi(A)=B$.  Put
\[ (A:B)_V:=|\det \phi|, \]
which is a positive real number independent of the choice of $\phi$. We call it
 the {\em Sinnott symbol} of $A$ to $B$. Note that 
\begin{enumerate}
\item[(1).] For lattices $A,\ B\subseteq V$, if $B\subseteq A$, 
then $(A:B)_V=\# (A/B)$. 
\item[(2).] Given lattices $A,\ B,\ C\subseteq V$, then $(A:B)(B:C)=(A:C)$.
\item[(3).] Let $f:V_1\rightarrow V_2$ be an isomorphism of vector spaces. Let $A$ and $B$
be lattices in $V_1$, then $(A:B)_{V_1}=(f(A):f(B))_{V_2}$. 
\end{enumerate}
For more results about the Sinnott symbol, see Sinnott~\cite{S1} and \cite{S2}.

Given a  finitely generated abelian group $A$, we denote the tensor product $A\otimes \R$ 
by $\R A$. Now given two finitely generated abelian groups $A$ and $B$, and an 
$\R$-linear isomorphism $\lambda:\R A\rightarrow \R B$. Choose free abelian 
subgroups $A'\subseteq A$ and $B'\subseteq B$ of finite index. Then 
$A'$ and $B'$ are of the same rank and hence isomorphic. Choose any isomorphism
$\phi:B'\rightarrow A'$, it can be naturally extended to an isomorphism
$\R\phi:\R B'\rightarrow\R A'$. make the evident identification $\R A'=\R A$
and $\R B'=\R B$. Now put
\begin{equation} 
\label{regdef}
\reg(A,B,\lambda):=\frac{|\det \R\phi\circ\lambda|\cdot \#B/B'}
{\#A/A'},\end{equation}
which is a positive real number independent of the choice of $A'$, $B'$ and $\phi$.
We call $\reg(A,B,\lambda)$ the {\em regulator} of $\lambda$ with respect to
$A$ and $B$. We often write it $\reg\lambda$ in abbreviation.
 
Here we calculate a few examples of the regulator:
\begin{exmp}
If both $A$ and $B$ are finite, then $ \reg(A,B,0)=\# B/\#A$.\end{exmp}

\begin{exmp} Let $f:A\rightarrow B$ be any homomorphism of finitely generated abelian
groups with finite kernel and cokernel, then $ \reg(A,B,\R f)$ is  exact
$\# \coker\ f/\# \ker f$.\end{exmp}

\begin{exmp} Let $A$, $B$ and $C$ be finitely generated abelian groups. Let $\lambda:
\R A\rightarrow \R B$ and  $\mu:\R B\rightarrow \R C$ be $\R$-linear isomorphisms. Then
$\reg \mu\circ\lambda=\reg \mu\cdot \reg \lambda$.\end{exmp}

\begin{exmp} Let $V$ be a finite dimensional $\R$-vector space. Let $A, B\subseteq V$ be
lattices. Let $\alpha:\R A\rightarrow V$ and $\beta:\R B\rightarrow V$ be the natural
isomorphisms induced by the inclusions $A \subseteq V$ and $ B\subseteq V$ 
respectively. Then $\reg (A, B, \beta^{-1}\circ\alpha)=(B:A)_V$.\end{exmp}

\subsection{The abstract index formula}
Consider  bounded complexes of finitely generated abelian groups
\[ (A,d_A): \cdots \rightarrow A^{i} \rightarrow A^{i+1}  \rightarrow  \cdots \]
and 
\[ (B,d_B): \cdots \rightarrow B^{i} \rightarrow B^{i+1}  \rightarrow  \cdots .\]
Given an isomorphism
\[ \lambda: \R A\longrightarrow \R B \] 
of bounded complexes of finitely dimensional vector spaces.  It naturally induces a map
\[ H^i(\lambda): H^i(\R A)\longrightarrow H^i(\R B) \]
for every degree $i$. Note that we also have $\R H^i(A)=H^i(\R A)$ and 
$\R H^i(B)=H^i(\R B)$. Then we have the following proposition:

\begin{prop} With the hypotheses above,  then
\begin{equation} \label{regco} \prod_i(\reg\ \lambda^i)^{(-1)^i}=
\prod_i(\reg\ H^i(\lambda))^{(-1)^i}. \end{equation}
\end{prop}
\begin{proof}
First we claim that there exist subcomplexes $A'\subseteq A$ and $B'\subseteq B$ satisfying 
the following conditions:
\begin{enumerate} \item[(1).]
 $A^{\prime i} $ and $B^{\prime i} $ are free abelian groups of the same rank as $A^i$ for 
all $i$;
\item[(2).]
 $H^i(A')$ and $H^i(B')$ are torsion free for all $i$;
\item[(3).]
 $A'$ and $B'$ are isomorphic complexes of abelian groups.
\item[(4).]
 The sequences 
\[ 0\rightarrow H^i(A')\rightarrow H^i(A)\rightarrow H^i(A/A')\rightarrow 0 \]
and 
\[ 0\rightarrow H^i(B')\rightarrow H^i(B)\rightarrow H^i(B/B')\rightarrow 0 \]
are exact for all $i$.
\end{enumerate}

\noindent This claim can be proved by induction. First since $A$ and $B$ are bounded 
complexes  of finite generated abelian groups, without loss of generality we suppose
\[ (A,d_A): \cdots 0\rightarrow A^{-n}\rightarrow\cdots\rightarrow A^{-1}
\rightarrow A^0\rightarrow 0\cdots \]
and
\[ (B,d_B): \cdots 0\rightarrow B^{-n}\rightarrow\cdots\rightarrow B^{-1}
\rightarrow B^0\rightarrow 0\cdots \]
Consider the subgroup $\im(d_A: A^{-1}\rightarrow A^0)$ of 
$A^0$. Let $r$ be the rank of $\text {im}\ A^{-1} $
 and let $\{e_1, \cdots, e_r\}$ be a maximal independent 
set in $\im\ A^{-1}$.  We can enlarge it into a maximal independent
set $E_0=\{e_1, \cdots, e_s\}$ of $A^0$. Set $A^{\prime 0}$ be the subgroup 
generated by $E_0$. Then $A^0/A^{\prime 0}$ is finite. Now consider 
the inverse image of $A^{\prime 0}$, it is a subgroup of $A^{-1}$. Moreover, 
it must have the same rank as $A^{-1}$. Since $\text{ker}(d_A: A^{-1}\rightarrow A^0)$ 
is contained in the inverse image of  $A^{\prime 0}$, so is $\im\ 
(d_A: A^{-2}\rightarrow A^{-1})$. 
Find $\{f_1, \cdots, f_s\} \subseteq A^{-1}$ such that $d_A(f_i)=e_i$. This set is
an independent set in the inverse image of  $A^{\prime 0}$ and has only trivial intersection
with  $\text{ker}(d_A: A^{-1}\rightarrow A^0)$.  
We select a maximal independent set in $\im(A^{-2}\rightarrow A^{-1})$, 
enlarge it to a maximal independent set in $\text{ker}(A^{-1}
\rightarrow A^0)$, together with $\{f_1, \cdots, f_s\} \subseteq A^{-1}$, we 
 get a maximal independent set $E_{-1}$ in the inverse image of 
$A^{\prime 0}$.  Set the free subgroup generated by $E_{-1}$ as 
$A^{\prime -1}$. Continuing this setup, we 
obtain a subcomplex $A'$ of $A$ such that $A^{\prime i}$ is free, $(A/A')^i $ is finite 
and $H^i(A')$ is torsion free.  

Similarly for the complex $B$, we can construct a subcomplex $B'$ of $B$ such that $
B^{\prime i}$ is free, $(B/B')^i $ is finite and $H^i(B')$ is torsion free.   Hence $A'$ and $B'$
satisfy conditions (1) and (2).  But (3) and (4) easily follow from (1) and (2). Hence we
proved the above claim.  Now choose an isomorphism $\phi: B'\rightarrow A'$ of complexes.
We have
\[  \begin{split} \prod_i(\reg\ \lambda^i)^{(-1)^i}=& \prod_i
\left (\frac{|\det \R\phi^i\circ \lambda^i|\cdot \#(B/B')^i}{\#(A/A')^i}\right )^{(-1)^i}\\
=& \prod_i\left (\frac{|\det \R H^i(\phi)\circ H^i(\lambda)|\cdot \#H^i(B/B')}
{\#H^i(A/A')}\right )^{(-1)^i}\\=& \prod_i(\reg\ H^i(\lambda))^{(-1)^i}.
\end{split} \]
Here we use the facts: (1). If $A$ is a complex of finite abelian group, then 
\[ \prod_i (\#H^i(A))^{(-1)^i} =\prod_i (\#A^i)^{(-1)^i}; \]
(2). If $V$ is a complex of $\R$-vector spaces, $\phi$ is an automorphism of $V$, then
\[ \prod_i \|\det \phi^i\|^{(-1)^i}=\prod_i \|\det H^i(\phi)\|^{(-1)^i}. \qed \]
\renewcommand{\qed}{} \end{proof}
Now Consider the following data:
\begin{itemize} \item A finite group $G$.
\item A bounded graded finitely generated left $\R[G]$-modules
\[ V=\bigoplus_i V^i\ \text{such that $V^i=0$ for $i>0$ and $i\ll 0$}, \]
equipped with two differential structures $d_1$ and $d_2$.
\item An  $\R[G]$-linear isomorphism $\phi$ between two cochain complexes
$(V,d_1)$ and $(V,d_2)$.
\item A lattice $L=\bigoplus_i L^i$ of $V$ which is $G$, $d_1$ and $d_2$-stable.
\item $H^i_{d_1}(L)=H^i_{d_2}(L)=0$ for all $i\neq 0$.
\item $H^0_{d_1}(L)$ and $H^0_{d_2}(L)$ are free abelian groups. 
\end{itemize}
Now for an arbitrary left ideal $\theta\subseteq \Z [G]$, by our assumption, we 
have the following trivial consequences:
\begin{itemize} \item $ H^i_{d_1} (V^{\theta})=H^i_{d_2} (V^{\theta})=0$ for all $i\neq 0$.
\item $L^{i \theta}$ is a lattice in $L^{i\theta}$ for all $i$.
\item $H^0_{d_2}(L)^{\theta}$ and $H^0_{d_2}(\phi L)^{\theta}$ are lattices in 
$H^0_{d_2}(V^{\theta})$.
\end{itemize}
By Proposition 2.1, we have(suggested by Anderson~\cite{ad1}
\begin{theo}[Abstract Index Formula] Under the above assumption, we have
\begin{equation} \label{aif} 
 (H^0_{d_2}(L)^{\theta}: H^0_{d_2}(\phi L)^{\theta})=\prod_i |\det(\phi^i\, |\, 
V^{i \theta})^{(-1)^i}|\cdot I(L,d_1;\theta)^{-1} \cdot I(L,d_2;\theta), \end{equation}
where for any complex of $\Z[G]$-modules $A$, we define 
\begin{equation} I(A;\theta):=\frac{\#\coker(H^0(A^{\theta})\rightarrow H^0(A)^{\theta})}
{\#\tor H^0(A^{\theta})\cdot \prod_{i\neq 0} \# H^i(A^{\theta})^{(-1)^i}} \end{equation}
if the above value is finite.
\end{theo}
\begin{proof} Consider the complexes $(L^{\theta},d_1)$ and $(L^{\theta},d_2)$ with the 
restriction map  $\phi: V^{\theta}\rightarrow V^{\theta}$.   Note that:

(1). $\reg(L^{i\, \theta}, L^{i\, \theta}, \phi^i)=|\det(\phi^i\, |\, V^{i \theta})|$ 
for all $i$.

(2). Since $ H^i_{d_1} (V^{\theta})=H^i_{d_2} (V^{\theta})=0$ for all $i\neq 0$, 
 $ H^i_{d_1} (L^{\theta})$ and $H^i_{d_2} (L^{\theta})$ are both finite and $H^i(\phi)=0$.
 We have $\reg(H^i_{d_1} (L^{\theta}), H^i_{d_2} (L^{\theta}), H^i(\phi))=
\# H^i_{d_2} (L^{\theta})/\# H^i_{d_1} (L^{\theta})$ for all $i\neq 0$.

(3). Now for $i=0$, consider the map $\alpha_j: H^0_{d_j} (L^{\theta})\rightarrow 
H^0_{d_j} (L)^{\theta}$.  We have $H^0(\phi)\circ \R \alpha_1=\R\alpha_2 \circ H^0(\phi)$.
Then 
\[ \begin{split} \reg(H^0_{d_1} (L^{\theta}),& H^0_{d_2} (L^{\theta}), H^0(\phi))\\
=&\reg(\alpha_1)\cdot \reg(\alpha_2)^{-1}
\cdot \reg(H^0_{d_1} (L)^{\theta}, H^0_{d_2} (L)^{\theta}, H^0(\phi)),
\end{split} \]
where 
\[ \reg(\alpha_j)=\frac{\# \coker(H^0_{d_j}  (L)^{\theta}\rightarrow H^0_{d_j} 
(L^{\theta}))} {\# \tor H^0_{d_j} (L^{\theta})} \]
and 
\[ \reg(H^0_{d_1} (L)^{\theta}, H^0_{d_2} (L)^{\theta}, H^0(\phi))
=(H^0_{d_2}(L)^{\theta}: H^0_{d_2}(\phi L)^{\theta}). \]

Now applying Formula~\eqref{regco} in Proposition 2.1 to the case $A=(L^{\theta},d_1)$, 
$B=(L^{\theta},d_2)$ and $\lambda=\phi$, we immediately get \eqref{aif}.
\end{proof}

\section{Theory of  spectral sequences }
\subsection{Basic theory of spectral sequences}
Let $G$ be a group and let $\Z[G]$ be the integral group ring of $G$. Let
$\theta$ be a left ideal of  $\Z[G]$. For any left module $M$, let $M^{\theta}$
be the subgroup of $M$ annihilated by $\theta$. Let
\[ (A,d):\cdots \rightarrow A^i\rightarrow A^{i+1}\rightarrow \cdots\]
be a complex of left $G$-modules. Assume
\begin{itemize} \item
 $A^i=0$ for $i>0$ and $i\ll 0$.
 \item
 $H^i(A)=0$ for $i\neq 0$.
\end{itemize}
\noindent Let $M=\Z[G]/\theta$, we have a projective resolution of $M$:
\[ (P,\partial):\cdots \rightarrow P_i\rightarrow \cdots\rightarrow P_{1}
\rightarrow P_0\rightarrow 0\]
Let $ K^{p,q}=\Hom_G(P_q, A^p)$, then we have a commutative
diagram:

\[
\begin{CD}
K^{p,q+1} @>{d\circ }>> K^{p+1,q+1}\\
@AA{\circ \partial}A @AA{\circ \partial}A  \\
K^{p,q} @>{d\circ }>> K^{p+1,q}
\end{CD}  \]
\vskip 0.4cm

\noindent    With  abuse of notations,
 we denote $d\circ$ by $d$ and $(-1)^p\ \circ \partial$
by $\delta$. Then we get a double complex $K^{\ast,\ast}=(K^{p,q}; d, 
\delta)$. The associate single complex is then defined by
\begin{equation} K^n=\bigoplus_{p+q=n} K^{p,q}, \ D=d+\delta. \end{equation}
Recall that we have two filtrations of the double complex $K^{\ast,\ast}$
\begin{equation}{'\Fil}^p K^{\ast,\ast}=\bigoplus_{p'\geq p} K^{p',q}, \end{equation}
and
\begin{equation} {''\Fil}^q K^{\ast,\ast}=\bigoplus_{q''\geq q} K^{p,q''}. \end{equation}
Now consider the following diagram:  
\vskip 0.5cm
\[
\begin{CD}
@. @AA{\delta}A  @AA{\delta}A @AA{\delta}A \\
 @>d>>   K^{p-1,q+1}@>{d}>> K^{p,q+1} @>{d}>> K^{p+1,q+1}@>d>>\\
  @.  @AA{\delta}A  @AA{\delta}A @AA{ \delta}A  \\
@>d>> K^{p-1,q}@>{d}>>           K^{p,q} @>{d}>> K^{p+1,q}@>d>>\\
 @.  @AA{\delta}A  @AA{\delta}A @AA{ \delta}A  \\
@>d>> K^{p-1,q-1}@>{d}>>           K^{p,q-1} @>{d}>> K^{p+1,q-1}@>d>>\\
@.  @AA{\delta}A  @AA{\delta}A @AA{\delta}A  
\end{CD}  \]
\vskip 0.5cm
\noindent We have
\begin{equation} H^q_{\delta}(K^{p,\ast})=\Ext_G^q(M,A^p), \end{equation}
and 
\begin{equation} H^p_{d}(K^{\ast,q})=
\begin{cases} 0,\ &\text{if}\ p\neq 0;\\ \Hom_G(P_q, H^0(A)),\ 
&\text{if}\ p= 0. \end{cases} \end{equation}
Therefore we can compute the $E_2$ terms of the related spectral sequences. 
 For the first one, 
\begin{equation} 'E_2^{p,q}=H^p(\Ext_G^q(M,A)); \end{equation}
for the second one, 
\begin{equation} ''E_2^{p,q}=
\begin{cases} 0,\ &\text{if}\ p\neq 0;\\ \Ext^q_G(M, H^0(A)),\ 
&\text{if}\ p= 0. \end{cases} \end{equation}
Since the second case collapses at $p=0$, we have
\begin{equation} H^i(K^{\ast})=\Ext_G^i(M,H^0(A)). \end{equation}
From now on we will focus only on the first
case. We omit the symbol $'$  from our notations. Then
\begin{equation} E_2^{p,q}=H^p(\Ext_G^q(M,A))\Rightarrow \Ext_G^{p+q}(M,H^0(A)).
 \end{equation}
Set $q=0$, then
\begin{equation} E_2^{p,0}=H^p(\Ext_G^0(M,A))=H^p(A^{\theta}).\end{equation}
Because $\Fil^1 K^{\ast}$ is trivial, we have
\[ E_{\infty}^{0,0}=\Fil^0 H^0(K^{\ast})=\im\ (H^0(\Fil^0 K^{\ast})
\rightarrow H^0(K^{\ast})). \]
Since $\Fil^0 K^{\ast}$ is nothing but the complex
\[ 0\rightarrow \Hom_G(P_0, A^0)\rightarrow \Hom_G(P_1, A^0)
\rightarrow \cdots \rightarrow \Hom_G(P_q, A^0)\rightarrow \cdots, \]
we have $H^0(\Fil^0 K^{\ast})=A^{0\, \theta}$ and
\[ E_{\infty}^{0,0}=\im\ (A^{0\, \theta}\rightarrow H^0(A)^{\theta}).\]
We show further it factors through $H^0(A^{\theta})$. First note that
$H^0(A^{\theta})=\coker(A^{-1\, \theta}\rightarrow A^{0\, \theta})$, 
therefore we only need to show that $A^{-1\, \theta}$ is contained in the 
boundary of $K^0$. This follows immediately from the diagram
\[
\begin{CD}
0 @>>>   A^{0\, \theta}@>>> K^{0,0} @>{\delta}>> K^{0,1}\\
  @.  @AA{d}A  @AA{d}A @AA{ d}A  \\
0 @>>>   A^{-1\, \theta}@>>> K^{-1,0} @>{\delta}>> K^{-1,1}
\end{CD}  \]
\noindent which is exact at the two rows. Combining the above arguments, 
we have
\begin{equation} E_{\infty}^{0,0}=\im\ (H^0(A^{\theta})\rightarrow H^0(A)^{\theta}).\end{equation}

\subsection{Application to the abstract index formula}
By the results obtained in the above subsection, we can express $I(A,\theta)$ in terms of
the order of $E_r$. We give here an important special case: 

\begin{prop} \label{lapse} If one has
\begin{equation} \label{epq}
\# \Ext_G^1(M,H^0(A))=\prod_q \# H^{1-q}(\Ext_G^q(M,A)), \end{equation}
then 
\begin{equation} \label{eppde}
I(A;\theta)=\prod_{\substack{p+q\leq 0\\q>0}} \#H^{p}(\Ext_G^q(M,A))^{(-1)^{p+q}}
=\prod_{\substack{p+q\leq 0\\q>0}} (\#E_2^{p,q})^{(-1)^{p+q}}. \end{equation}
\end{prop}
\begin{proof}
First note that the given identity \eqref{epq} is nothing but
\[ \prod_q \#E_{\infty}^{1-q,q}=\prod_q \#E_{2}^{1-q,q}. \]
Since for the spectral sequence, $H^{\ast}(E_r)=E_{r+1}$, we always have
\[\#E_{2}^{p,q}\geq \#E_{3}^{p,q}\geq\cdots \geq \#E_{\infty}^{p,q}. \]
Hence 
\[ \#E_{2}^{1-q,q}= \#E_{3}^{1-q,q}=\cdots = \#E_{\infty}^{1-q,q}, \]
which means that for  $r\geq 2$,
\[ \im(d_r: E_r^{1-q-r, q+r-1}\rightarrow E_r^{1-q, q})=\im(d_r: 
E_r^{1-q, q}\rightarrow E_r^{1-q+r, q-r+1})=0.\]
Therefore we have a shorter complex:
\[ \cdots\rightarrow E_r^{1-q-2r, q+2r-2}\rightarrow E_r^{1-q-r, q+r-1}
\rightarrow 0. \]
Now we set to prove the following fact:
\begin{equation} \label{eppp}
\prod_{\substack{p+q\leq 0\\(p,q)\neq (0,0)}} (\#E_r^{p,q})^{(-1)^{p+q}}
\cdot \# \tor E_r^{0,0}=\text{Constant}. \end{equation}
Observe that the set $\{ E_r^{p,q}: p+q\leq 0,\ q\geq 0\}$, the only term not 
finite is $ E_r^{0,0}$. If we substitute it by its torsion, we still get 
a group of complexes composed of finite abelian groups and with differential 
$d_r$. The cohomology groups are $E_{r+1}^{p,q}$(or $\tor E_{r+1}^{0,0}$).
By the invariance of Euler characteristic under cohomology, \eqref{eppp} is proved.
Note that $E_{\infty}^{0,0}$ is free and
\[ \prod_{\substack{p+q\leq 0\\(p,q)\neq (0,0)}} 
(\#E_{\infty}^{p,q})^{(-1)^{p+q}} 
=\#\coker(H^0(A^{\theta})\rightarrow H^0(A)^{\theta}). \]
The formula~\eqref{eppde} now follows immediately.
\end{proof}

\section{The universal distribution and predistribution}
\subsection{Definitions and basic properties} Let $\A$ be the free abelian group generated by the
symbols $[a]$ with $a\in  \Q/\Z$. We call the elements which are linear combinations of 
 $[a]-\sum_{nb=a} [b]$ distribution relations in $\A$ and 
the elements which are linear combinations of $\sum_{nb=a} [b]$
predistribution relations in $\A$.  Let  $U$ be the quotient group of $\A$ modulo
the distribution relations and let $\co$ the quotient group of $\A$ modulo the predistribution 
relations. We call $U$ and $\co$ the (rank 1) \em{universal distribution} and  the (rank 1) 
\em{universal predistribution}
respectively. Now for the subgroup $\A_m=<[a]: a\in \frac{1}{m}\Z/\Z>$ of $\A$, put
\begin{align}  
 U_m&=\A_m/<[a]-\sum_{nb=a}[b], n|m,\ a\in \frac{n}{m}\Z/\Z>,\notag  \\
\intertext{ and}
 \co_m&=\A_m/<\sum_{nb=a}[b], n|m,\ a\in \frac{n}{m}\Z/\Z>. \notag \end{align}
We call $U_m$ and $\co_m$ the universal distribution and the universal 
predistribution of level $m$ respectively.

In \cite{S1}, Sinnott introduced an $\R[G]$-module $U_m$ and used it to compute the index
of the Stickelberger ideal and  the circular units.  Kubert \cite{Kubert1} then proved that
Sinnott's module are actually isomorphic to the one we defined above. We have 

\begin{prop}
\[ \begin{split}
(1).\ U_m\ \cong\ & U_{\text{Sinnott}};\\
(2).\ \co_m\ \cong\ & \co_{K_m}. \end{split} \] 
\end{prop}
\begin{proof} (1). See Kubert \cite{Kubert1}.

(2). Define
\[ \begin{split} e_m:\ \A_m\ &\longrightarrow\ \co_{K_m}\\
\sum n_i[a_i]&\longmapsto \sum n_i\exp(2\pi i a_i) \end{split}\]
It is routine to check that $e_m$ is actually an isomorphism. Since we don't need this fact
in the latter context, we omit it here.
\end{proof}
\subsection{The connecting map $\phi_m$} Let
\[ \begin{split} 
\phi_m: \R\otimes\A_m&\longrightarrow \R\otimes \A_m\\
[x]&\longmapsto \sum_{n\mid m^{\infty}} \frac{[nx]}{n}.
\end{split} \]
Then $\phi_m$ is an automorphism of $\R$-vector space $\R\A_m$,  the inverse map is 
given by 
\[ \phi_m^{-1}:[x]\longmapsto \sum_{n\mid m^{\infty}} \frac{\mu(n)[nx]}{n}, \]
where 
\[ \mu(n)=\begin{cases} (-1)^i,\ &\text{if  $n$ is a product of $i$ distinct
prime numbers};\\0,\ &\text{otherwise}. \end{cases} \]
is the M\"obius function. Now if we enlarge the definition of distribution and predistribution
relations to $\R\A_m$, then we have
\begin{prop}  \label{connectionth} 
$\phi_m$ maps distribution relations to predistribution relations. In other words,
$\phi_m$ induces an isomorphism from $\R U_m$ to $\R \co_m$.
\end{prop}
\begin{proof} By straightforward calculation. \end{proof} 
\begin{note} From now on we denote by $\varphi_m$  the above induced map. \end{note}

\section{The cochain complexes $(L_m, d_{1m})$ and  $(L_m, d_{2m})$}
\subsection{Set up}In this section and sequel, we  fix the 
following notations:
\begin{itemize} 
\item \( K_m=\Q(\zeta_m),\ G_m=\text{Gal}(\Q(\zeta_m)/\Q)=(\Z/m\Z)^{\times}; \)
\item $c=\sigma_{-1}$ is the complex conjugation in $G_m$, $\theta=1+c$, $J=\{1,c\}$.
\item \( L_m= <[x,g]: g|m,\ g\ \text{square free},\ x\in \frac{g}{m}\Z/\Z>; \)
\item \( L_{m,g}=<[x,g]: x\in \frac{g}{m}\Z/\Z>\) for a fixed square free factor $g$; 
\item \( L_m^i= \bigoplus L_{m,g}\) for all square free $g|m$ such that $i=-\#
\text{Supp}\ g$;
\item \( V_m=\R\otimes L_m,\ V_{m,g} =\R\otimes L_{m,g},\ V^i_m= \R\otimes L^i_m . \)
\end{itemize}
\noindent For any square free positive integer $g$, suppose that $g=p_i\cdots 
p_r$,  $p_1<\cdots <p_r$, is the prime factorization of $g$. Put
\[ \epsilon(g,p)=\begin{cases} (-1)^i,\ & \text{if}\ p=p_i;\\
0,\ & \text{otherwise}. \end{cases} \]
Now we define 
\begin{equation} d_{1m}: L^i_m\rightarrow L^{i+1}_m, [x,g] \mapsto \sum_{i=1}^r 
 \epsilon(g, p_i)   ([x, g/p_i]-\sum_{p_i y=x}[y, g/p_i]). \end{equation}
and
\begin{equation} d_{2m}: L^i_m\rightarrow L^{i+1}_m, [x,g] \mapsto \sum_{i=1}^r \epsilon
(g, p_i) (-\sum_{p_i y=x}[y, g/p_i]), \end{equation}
By straightforward calculation, we have $d^2_{1m}=d^2_{2m} =0$. 
Therefore $V_m$ is equipped with two cochain complexes structure,  we write them 
$(V_m, d_{1m})$ and $(V_m, d_{2m})$ respectively..
In the next section, we are going to study the cohomology groups. 
\subsection{Connecting map again}  
In this subsection,  we define a connecting map $\phi_m$ between $(V_m, d_{1m})$ and 
$(V_m, d_{2m})$, generalizing the one defined in $\S 4.2$. We put 
\[ \begin{split} 
\phi_m: V_m&\longrightarrow V_m\\
[x,g]&\longmapsto \sum_{\substack{n\mid m^{\infty}\\(n,g)=1}} \frac{[nx,g]}{n}.
\end{split} \]
Then $\phi_m$ is an automorphism of the vector space $V_m$. Furthermore,  the inverse 
map of $\phi_m$ is given by
\[ \begin{split} 
\phi_m^{-1}: V_m&\longrightarrow V_m\\
[x,g]&\longmapsto \sum_{\substack{n\mid m^{\infty}\\(n,g)=1}} \frac{\mu(n)[nx,g]}{n}.
\end{split} \]
The following proposition establishes the connection between  $(V_m, d_{1m})$ and 
$(V_m, d_{2m})$.
\begin{prop} $ \phi_m$ is an isomorphism from cochain complex $(V_m, d_{1m})$
to cochain complex  $(V_m, d_{2m})$, i.e. , 
\[ d_{2m}\phi_m=\phi_m d_{1m}. \]
\end{prop}
\begin{proof} By direct calculation.
\end{proof}
\begin{rem} Under the apparent isomorphism from $V^{0}_m$ to $\R\otimes \A_m$, 
we can see that the connecting map $\phi$ defined in $\S 3.2$ is the same map $\phi$ defined 
on $V^{0}_m$. Later we will see that $\varphi=H^0(\phi)$. 
\end{rem}
Now we try to calculate the determinant of $\phi_m$.
We have 

\begin{prop} \label{deter}
 \begin{equation} \label{amsmam} \prod_{i}\det(\phi_m: V^i_m)^{(-1)^i}=
\prod_{p\mid m}\; \prod_{\chi\in \hat{G}_m}(1-\chi(p)p^{-1})^{-1}. \end{equation}
\end{prop}
\begin{proof} First notice that $V_{m,g}$ is invariant under $\phi_m$. 
Moreover, let $h=m/g$, for any $f\mid h$, define
\[ V^f_{m,g}=\R\, \otimes <[x,g]: fx=0>, \]
then clearly $V^f_{m,g}$ is invariant under $\phi_m$. By definition, 
we have $ V^h_{m,g}=V_{m,g}$. Put 
\[ V^{(f)}_{m,g} = V^{f}_{m,g}/\sum_{p\mid f}V^{f/p}_{m,g}, \]
We can see that $ V^{(f)}_{m,g}$ is a real vector space with a basis 
$\{ [\frac{a}{f},g]:(a,f)=1\}$. Furthermore $ V^{(f)}_{m,g}$ has a natural
$\R[G_f]$-module structure. Actually it is a free $\R[G_f]$-module of
rank $1$.  $\phi_m$ induces an automorphism in $ V^{(f)}_{m,g}$:
\[ \begin{split}
\phi_m:  V^{(f)}_{m,g} &\longrightarrow V^{(f)}_{m,g}\\
[x,g]&\longmapsto \sum_{\substack{n\mid m^{\infty}\\(n,fg)=1}} \frac{[nx,g]}{n}.
\end{split} \]
We calculate its determinant first. Let
\[ S_{f,g}=\text{Supp}\ m- \text{Supp}\ f\cup\text{Supp}\ g. \]
For $p\in S_{f,g}$, define
\[ \begin{split}
\tau_p: V^{(f)}_{m,g} &\longrightarrow V^{(f)}_{m,g}\\
[x,g]&\longmapsto \sum_{{n\mid p^{\infty}}} \frac{[nx,g]}{n}.
\end{split} \]
Note that $ \tau_{p_i}\circ\tau_{p_j}= \tau_{p_j}\circ\tau_{p_i}$ and 
\[ \phi_m|_{ V^{(f)}_{m,g} }=\tau_{p_1}\circ\cdots \circ \tau_{p_s} \]
where $p_i\in S_{f,g}$. Then we have
\[ \det(\phi_m:\ {V^{(f)}_{m,g}})=\prod_{p\in S_{f,g}}\det \tau_p. \]
For any $p\in S_{f,g}$, let $c_{p,f}$ be the smallest number satisfying 
$p^{c_{p,f}}\equiv 1\pmod f$. Since the map $\tau_p$ can be regarded
as the left multiplication by the group ring element $\sum_i \frac{\sigma^i_p}
{p^i}$ in $\R[G_f]$, then by \cite{S2} Lemma 1.2(b), we have
\[ \det \tau_p=\prod_{\chi\in \hat{G}_f}\chi(\sum_i \frac{\sigma^i_p}{p^i})=
(1-p^{-c_{p,f}})^{-\varphi(f)/c_{p,f}}:=a_{p,f}, \]
and 
\[ \det(\phi_m:\ V^{(f)}_{m,g})=\prod_{p\in S_{f,g}}
a_{p,f}. \]
Now by the Inclusion-Exclusion Principle, we have
\[ \det(\phi_m:\sum_{p\mid f}V^{f/p}_{m,g})=\prod_{f'\mid f, f'\neq 1}
\det(\phi_m: V^{f/f'}_{m,g})^{-\mu(f')}. \]
Hence 
\[ 
\prod_{p\in S_{f,g}} a_{p,f}
=\prod_{f'\mid f} \det(\phi_m: V^{f/f'}_{m,g})^{\mu(f')}. \] 
By the M\"obius inverse formula, 
\[ \det(\phi_m:V_{m,g})=\prod_{f\mid \frac{m}{g}} 
\prod_{p\in S_{f,g}} a_{p,f}  \]
Therefore we have
\[ \prod_{i} \det(\phi_m: V^i_m)^{(-1)^i}
=\prod_{g\mid m} \Bigl (\prod_{f\mid \frac{m}{g}} 
\prod_{p\in S_{f,g}} a_{p,f}\Bigr )^{\mu(g)}. \]
Now let's look at the right hand side of the above identity. The exponent
of $a_{p,f}$ is
\[ 
\sum_{{g\mid \frac{m}{f},\ (p,g)=1}} \mu(g)=\sum_{g\mid \frac{m}
{fp^{\alpha}}} \mu(g)
=\begin{cases}1,\ & \text{if}\  \frac{m}{fp^{\alpha}}=1;\\
0,\ & \text{otherwise}. \end{cases}  \]
here $p^{\alpha}\|\, m$. Write $m=m_p\cdot p^{\alpha}$, then 
\[ \prod_{i}\det(\phi_m: V^i_m)^{(-1)^i}=\prod_{p\mid m} a_{p, m_p}
=\prod_{\mathfrak P\cap \Z=p\mid m}\Bigl (1-\frac{1}{N\mathfrak P}
\Bigr )^{-1}, \]
which is exact the right hand side of  the identity~\eqref{amsmam}. \end{proof}
Now let $V^{\theta}_m=\{x\in V_m: \theta\cdot x=0\}$. 
$V^{\theta}_m$ has a basis consisting of $\{[x,g]-[-x,g]: 0< x<1/2\}\subseteq V_m$.
Denote by $\phi^{\theta}_m$ the restriction of $\phi_m$ on $V^{\theta}_m$. Then 
$\phi^{\theta}_m$ is an automorphism of $V^{\theta}_m$. We have
\begin{prop} \label{deter1}
 \begin{equation} \prod_{i}\det(\phi^{\theta}_m: V^{i\theta}_m)^{(-1)^i}=
\prod_{\chi\ odd}\prod_{p\mid m} (1-\chi(p) p^{-1})^{-1}. \end{equation}
\end{prop}
\begin{proof} The proof is similar to the proof of Proposition~\ref{deter}. Note that
$ V^{(f)\theta}_{m,g}$ is a real vector space with a basis 
$\{ [\frac{a}{f},g]-[-\frac{a}{f},g]:(a,f)=1, 0<a<f/2\}$. On the quotient space $ V^{(f)\theta}_{m,g}$,
\[ \phi^{\theta}_m: [x,g]-[-x,g]\longmapsto \sum_{\substack{n\mid m^{\infty}\\(n,fg)=1}} 
\frac{[nx,g]-[-nx,g]}{n}.  \]
Now the restriction of $\tau_p$ on  $V^{(f)\theta}_{m,g}$ is
\[ \tau^{\theta}_p: 
[x,g]-[-x,g]\longmapsto \sum_{{n\mid p^{\infty}}} \frac{[nx,g]-[-nx,g]}{n}. \]
We still have 
\[ \phi^{\theta}_m|_{ V^{(f)\theta}_{m,g} }=\tau^{\theta}_{p_1}\circ\cdots \circ 
\tau^{\theta}_{p_s} \]
where $p_i\in S_{f,g}$. Similar to the calculation of $\det\tau_p$ in Proposition~\ref{deter}, 
we have
\[ \det \tau^{\theta}_p:=b_{p,f}=\begin{cases}
(1-p^{-c_{p,f}})^{-\varphi(f)/2c_{p,f}},\ & if \ c_{p,f}\ odd;\\
(1+p^{-c_{p,f}/2})^{-\varphi(f)/c_{p,f}}, \ & if \ c_{p,f}\ even. \end{cases} \]
We have
\[ 
\prod_{i}\det(\phi^{\theta}_m: V^{i\theta}_m)^{(-1)^i}=\prod_{p\mid m} b_{p, m_p}
=\prod_{\chi\ odd}\, \prod_{p\mid m} (1-\chi(p) p^{-1})^{-1}.  \qed \]
\renewcommand{\qed}{} \end{proof}

\section{ Computation of  $H^{\ast}(L_m, d_{1m})$ and $H^{\ast}(L_m, d_{2m})$}
This section is dedicated to the computation of the cohomology groups 
$H^{\ast}(L_m, d_{1m})$ and $H^{\ast}(L_m, d_{2m})$. We first introduce module structures on
$\A$ and $L=\cup L_m$. Using this structure, we find new bases for $\A$ and $L$,
which is applied to study the cohomology groups.

Let $\Lambda= \Z [X_2, X_3,\cdots,X_p,\cdots] $
be the polynomial ring generated by indeterminants $X_p$ for all prime number 
$p$. For every positive integer $n=\prod p^{n_p}$, put
\[ X_n=\prod {X_p}^{n_p}, \qquad  Y_n=\prod (1-{X_p})^{n_p}, \]
then the set $\{X_n, n\in \N\}$ is a $\Z$-basis of $\Lambda$, so is $\{Y_n, 
n\in \N\}$. Now $\A$ and $L$ are equipped with $\Lambda$-module structures  by 
the following rules:
\[ X_n[a] =\sum_{nb=a}[b],\ X_n[a,g] =\sum_{nb=a}[b,g]. \]
Put
\[ d_1[a,g]:= \sum_{p\mid g} \epsilon(g,p)Y_p[a, g/p], \]
and 
\[ d_2[a,g]:= \sum_{p\mid g} \epsilon({g,p})\cdot (-X_p)[a, g/p]. \]
Then $d_1^2=d_2^2=0$ and $L$ is equipped with a cochain complex 
structure by $d_1$ or $d_2$. Furthermore 
\[ d_1|_{L_m}=d_{1m}, \qquad  d_2|_{L_m}=d_{2m}. \]

For any $a\in\Q/\Z$, we can uniquely write 
\[ a\equiv \sum_p \sum_v \frac{a_{pv}}{p^v}\pmod \Z, \]
where $0\leq a_{pv}<p$ for each pair of any prime number $p$ and any positive
 integer $v$. Note that $a_{pv}=0$ for all but finite any $\{p,v\}$.
For each nonnegative integer $k$ we define $\mathcal R_k$ to be
 the set of $a\in \Q/\Z$ such that there exist at most $k$ prime numbers $p$
such that $a_{p1}=p-1$. In particular, $\mathcal R_0$ is  the set of $a\in 
\Q/\Z$ such that such that $a_{p1}\neq p-1$ for all prime numbers $p$.
\begin{prop} \label{e}
(1).  For each positive integer $m$, the collection
\[ \{ X_n[a]: n\mid m, \ n\in \N, a\in \mathcal R_0\cap \frac{n}{m}\Z/\Z\} \]
constitutes a basis for the free abelian group $\A_m$.

(2). The collection $ \{ X_n[a]: \ n\in \N, a\in \mathcal R_0\} $
constitutes a basis for the free abelian group $\A$.

(3). As a $\Lambda$-module $\A$ is free with a $\Lambda$-basis $\{ [a]:  
a\in \mathcal R_0\}$.

(4). In (1) and (2), if we change $X_n$ by $Y_n$, the related results are 
still true.
\end{prop}
\begin{proof}
First note that $(1)\Rightarrow (2)\Rightarrow (3)$. For (1), since 
\[ |\mathcal R_0\cap \frac{1}{m}\Z/\Z|=\varphi(m)=|(\Z/m\Z)^{\times}|, \]
it suffices to show that the given collection generates $\A_m$. This can be easily
deduced by induction to $\mathcal R_k$, with the fact that for any $a\in \A_m$,
\[ [a]=-\sum_{i=1}^{p-1}[a-\frac{i}{p}]+ X_p[pa]. \]
For (4), note that the identity
\[ X_n- (-1)^{\sum n_i} Y_n =\sum_{l\mid n} c_{nl} X_l, \]
holds for any $n=\prod p_i^{n_i}$ and integer constants $c_{nl}$, therefore (4) follows
 immediately from (1) and (2).
\end{proof}
By Proposition~\ref{e}, we have
\begin{prop} \label{f}
(1).  For each positive integer $m$, the collection
\[ \{ X_n[a,g]: ng\mid m, \ n\in \N, \ g\ squarefree , a\in \mathcal R_0\cap 
\frac{ng}{m}\Z/\Z\} \]
constitutes a basis for the free abelian group $L_m$.

(2). The collection $ \{ X_n[a,g]: \ n\in \N, \ g\ squarefree,\ a\in 
\mathcal R_0\} $ constitutes a basis for the free abelian group $L$.

(3). In (1) and (2), if we change $X_n$ by $Y_n$, the related results are 
still true.
\end{prop}

With the help of Proposition~\ref{f}, we can compute the cohomology groups of 
$(L_m, d_{1m})$ and $(L_m, d_{2m})$.  
\begin{theo} \label{g} (1). The complex $(L, d_1)$ and $(L, d_2)$ are acyclic 
in negative degree, moreover, $H^0(L, d_1)$ is the universal distribution $U$
and $H^0(L, d_2)$ is the universal predistribution $\co$.

(2). The complex $(L_m, d_{1m})$ and $(L_m, d_{2m})$ are acyclic 
in negative degree, moreover, $H^0(L_m, d_{1m})$ is $U_m$
and $H^0(L_m, d_{2m})$ is $\co_m$. In both cases, the
natural map $H^0(L_m)\rightarrow H^0(L)$ is injective.
\end{theo}
\begin{proof} For each prime number $p$, we define operators $d_2^p$,
$T_2^p$ and $\pi_2^p$ on $L$ by the rules:
\[ \begin{split} &d_2^p X_n[a,g]=
\begin{cases} -\epsilon(g,p) X_{np}[a,g/p], \ &\text{if} \ p\mid g,\\ 0,\ & \text{otherwise}.
\end{cases} \\
& T_2^p X_n[a,g]=
\begin{cases} -\epsilon(gp,p) X_{n/p}[a,gp], \ &\text{if} \ (p,g)=1\ \&\ p\mid n,
\\ 0,\ & \text{otherwise}.\end{cases} \\
& \pi_2^p X_n[a,g]=
\begin{cases} X_{n}[a,g], \ &\text{if} \ (p,ng)=1,\\ 0,\ & \text{otherwise}.
\end{cases} \end{split} \]
where $X_n[a,g]$ runs through the basis given in Theorem ~\ref{f}(2).  It is easy to check:
\[ d_2^p T_2^q+T_2^q d_2^p =\delta_{pq}(1-\pi_2^p) \]
where $\delta_{pq}$ is the Kronecker symbol,  and 
\[ \pi_2^{p\, 2}=\pi_2^p,\qquad  \pi_2^p\pi_2^q=\pi_2^q\pi_2^p. \]

By Proposition~\ref{f}(1), for any fixed positive integer $m$, we can check that 
\[ d_2^p L_m, \ T_2^p L_m,\ (1-\pi_2^p) L_m \subseteq \begin{cases} 
L_m,\ & \text{if}\ p\mid m,\\ 0,\ & \text{if}\ (p,m)=1. \end{cases} \] 
Now for a fixed positive integer $m$, we define $d_{2m}$, $T_{2m}$ and $\pi_{2m}$ on 
$L_m$ by the rules:
\[ \begin{split} &d_{2m} :=\sum_{p\mid m} d_2^p, \\
& T_{2m}:=\sum_{p\mid m} (\sum_{q<p} \pi_2^q)T_2^p, \\
& \pi_{2m}:=\prod_{p\mid m} \pi_2^p, \end{split} \]
here we abuse the notations $d_2^p$, $T_2^p$ and $\pi_2^p$ with their restrictions on
$L_m$. It is easy to check the definition $d_{1m}$ here coincides the one we defined 
in $\S 4$. We can also check that
\[ \pi_{2m} X_n[a,g] =\begin{cases} [a,1],\ & \text{if} \ n=1\ \&\  g=1;\\ 0, \ &
\text{otherwise}.\end{cases} \]
for any elements $X_n[a,g] $ of the basis given in Theorem~\ref{f}(1). Now we have
\[ \begin{split} d_{2m}T_{2m} +T_{2m} d_{2m}
=& (\sum_{\ell\mid m} d_2^{\ell} )(\sum_{p \mid m}(\sum_{q<p} \pi_2^{q}) T_2^{p} )
+(\sum_{p \mid m}(\sum_{q<p} \pi_1^{q}) T_2^{p} )(\sum_{\ell\mid m} d_2^{\ell} )\\
=&\sum_{p}\sum_{\ell}(\sum_{q<p} \pi_2^{q})(d_2^{\ell}T_2^{p} +T_2^{p}d_2^{\ell})\\
=&\sum_{p}(\sum_{q<p} \pi_2^{q})(1-\pi_2^p)\\
=&\sum_{p}(\sum_{q<p} \pi_2^{q}-\sum_{q\leq p} \pi_2^{q})\\
=&1-\pi_{2m}. \end{split} \]
By the above argument, the cochain map $\text{id}: (L_m,d_{2m})\rightarrow (L_m,d_{2m})$ is 
homotopic to the cochain map $\pi_{2m}$.  But in the negative degree,  
$pi_{2m}$ is nothing but the zero map. Therefore we showed
 that $(L_m, d_{2m}) $ is acyclic in negative degrees. Now since for the map $d_{2m}$(resp.
$T_{2m}$, $\pi_{2m}$),  $d_{2m}|_{L_m\cap L_{m'}}=d_{2m'}|_{L_m\cap L_{m'}}$(resp.
$T_{2m'}$, $\pi_{2m'}$),  there exists a unique operator $d_{2}$(resp.
$T_{2}$, $\pi_{2}$) with restriction at $L_m$ the  operator $d_{2m}$
(resp.$T_{2m}$, $\pi_{2m}$). $\pi_{2}$ is a cochain map homotopic to the 
identity map of the cochain complex $(L,d_2)$ and vanishes at negative degrees. Therefore 
$(L, d_2)$ is acyclic at negative degrees. Now for $n=0$, the cohomology groups $H^0(L,d_2)$
and $H^0(L_m,d_{2m})$  easily follow from the definitions of $d_2$ and $d_{2m}$ .

Now by a parallel argument to $d_1$ and $Y_n$,  we construct $T_1$, $\pi_1$ and
$T_{1m}$, $\pi_{1m}$ respectively.  the remaining assertions  follow immediately
\end{proof}
\begin{rem} 1. The proof here is given by Anderson in the preprint version of \cite{ad2}.

2. In the higher rank case, we also have similar result  by applying essentially the same 
trick. 
\end{rem}
\section{More spectral sequences}
In the following sections we are going to use the spectral sequence method 
to attack the cochain complexes introduced in  $\S 4$.  First recall:

\begin{itemize}
\item $J=\Z/2\Z=\{1,c\}\subseteq G_m$,  $\theta=1+c\in \Z[J]$;  
\item  $L_m=<[a,g]:\, g\mid m,\, a\in\frac{g}{m}\Z/\Z>$,   $d=d_1$ or $d_2$; 
\item $H^0_{d_1}(L_m)=U_m$,  $H^0_{d_2}(L_m)=\co_m$;
\item $r=\# \text{Supp}\ m$, $g\mid m$, $g$ square free, $p=-\# \text{Supp}\ g$.
\end{itemize}
Now let $M=\coker(\Z[J]\stackrel{1+c} \rightarrow \Z[J])$. Consider two chain complexes 
\[ (P,\partial): \cdots \stackrel{\partial_{q+1}}\longrightarrow\Z[J]_{q+1}\stackrel{\partial_q}
 \longrightarrow \Z[J]_q\stackrel{\partial_{q-1}}\longrightarrow
\cdots\stackrel{\partial_0}\longrightarrow \Z[J]_0 \longrightarrow 0 \]
and
\[ (F,\partial): \cdots \stackrel{\partial_{q+1}}\longrightarrow\Z[J]_{q+1}\stackrel{\partial_q} 
\longrightarrow \Z[J]_q \stackrel{\partial_{q-1}}\longrightarrow\cdots\stackrel{\partial_0}
\longrightarrow \Z[J]_0 \stackrel{\partial_{-1}} \longrightarrow \cdots \]
where $\Z[J]_q=\Z[J]$ and $\partial_{q}=1+(-1)^q\cdot c$. It is clear that
$P$ is a projective resolution of $M$ and $F$ is an exact sequence. Regard the cochain
complex $(L_m,d)$ as $A$ in $\S 2$, we can construct two double complexes by
\[ K^{p,q}=\begin{cases} \Hom_G(\Z[J]_q, L_m^p):= (L_m^p,q),\ & \text{if}\  q\geq 0
\\ 0,\ & \text{if}\ q< 0 \end{cases}\]
and 
\[ F^{p,q}= (L_m^p,q) \]
where the induced differentials $\delta_q: (x,q)\mapsto ((-1)^p(1+(-1)^q c)x,q+1)$ and 
$d: (x,q) \mapsto (d(x), q)$. The two filtrations for the double complex
$K^{\ast,\ast}$ are
\[ {'\Fil}^p K^{\ast,\ast}=\bigoplus_{p'\geq p} K^{p',q} \]
and
\[ {''\Fil}^q K^{\ast,\ast}=\bigoplus_{q''\geq q} K^{p,q''}. \]
From results in $\S 2$, we have 
\[ ''E_2^{p,q}=
\begin{cases} 0,\ &\text{if}\ p\neq 0;\\ \Ext^q_J(M, H_d^0(L_m)),\ 
&\text{if}\ p= 0. \end{cases} \]
Using the projective resolution $P$ of $M$,  the  cohomology groups of the total complex 
$K^{\ast}$ are 
\[ H^n(K^{\ast})=\Ext_J^n(M,H_d^0(L_m))=\begin{cases} 
 H^1(J, H_d^0(L_m)), \ &\text{if $n$ odd},\ n>0;\\
H^2(J, H_d^0(L_m)), \ &\text{if $n$ even},\ n>0;\\
H_d^0(L_m)^{1+c},\ &\text{if }\ n=0;\\
0,\ &\text{if }\ n<0. \end{cases} \]
Now consider the first filtration,  we have 
\[ 'E_1^{p,q}=H_{\delta}^q(K^{p,\ast})=\begin{cases} 
 H^1(J, L_m^p), \ &\text{if $q$ odd},\ q>0;\\
H^2(J, L_m^p), \ &\text{if $q$ even},\ q>0;\\
(L_m^p)^{1+c},\ &\text{if }\ q=0;\\
0.\ &\text{if }\ q<0. \end{cases} \]
Similarly for the double complex $F^{\ast,\ast}$, 
\[ H^q(F^{\ast})=\hat{H}^q(J,H_d^0(L_m))=\begin{cases} 
 H^1(J, H_d^0(L_m)), \ &\text{if $q$ odd};\\
H^2(J, H_d^0(L_m)), \ &\text{if $q$ even}. \end{cases} \]
and
\[ 'E_{1,F}^{p,q}=H_{\delta}^q(F^{p,\ast})=\hat{H}^q(J,L_m^p)=\begin{cases} 
 H^1(J, L_m^p), \ &\text{if $q$ odd};\\
H^2(J, L_m^p), \ &\text{if $q$ even}. \end{cases} \]
Now for the Galois cohomology $\hat{H}^q(J, L_{m}^p)$, first since(Note that 
$m$ $\not \equiv 2$ mod $4$ by our assumption)  
\[  L_{m,g}=\begin{cases}
\bigoplus_{2a\neq 0}(\Z[a,g]\bigoplus \Z[-a,g])\bigoplus \Z[0,g]\bigoplus 
\Z[\frac{1}{2},g],\ & \text{if $m$ even};\\
\bigoplus_{2a\neq 0}(\Z[a,g]\bigoplus \Z[-a,g])\bigoplus \Z[0,g],\ & \text
{if $m$ odd}. \end{cases} \]
Then 
\[ \hat{H}^q(J, L_{m,g})=\begin{cases} 
(\Z/2\Z)^2, \ & \text{if $q$ odd, $m$ even};\\
(\Z/2\Z), \ & \text{if $q$ odd, $m$ odd};\\
0, \ & \text{if $q$  even}.\end{cases} \]
Hence 
\[ \hat{H}^q(J, L_{m}^p)=\begin{cases} 
(\Z/2\Z)^{2\binom{r}{-p}}, \ & \text{if $q$ odd, $m$ even};\\
(\Z/2\Z)^{\binom{r}{-p}},\ & \text{if $q$ odd, $m$ odd};\\
0, \ & \text{if $q$  even}.\end{cases} \]
Denote by $X_g$ the cocycle represented by $[0,g]$ and 
 by $Y_g$ the cocycle represented by $[1/2,g]$, then for $q>0$,
\[  'E_1^{p,q}=\hat{H}^q(J, L_{m}^p)==\begin{cases} 
\bigoplus_g (<X_g>\bigoplus <Y_g>), \ & \text{if $q$ odd, $m$ even};\\
\bigoplus_g  <X_g>, \ & \text{if $q$ odd, $m$ odd};\\
0, \ & \text{if $q$  even}.\end{cases} \]
We consider $ 'E_1^{p,q}$ as a finite dimensional $\Z/2\Z$ vector space.
Immediately we have $ 'E_2^{p,q}=0$ for $q$ even. Now for $q$ odd,
if $m$ is odd, the induced differential $d^1$ is
\[  X_g\stackrel{d_1^1}  \longmapsto 0,  \qquad  X_g\stackrel{d_2^1}
 \longmapsto\sum_{i=1}^{-p}X_{g/p_i}; \]
if $m$ is even,  the induced differential $d^1$ is
\[  X_g\stackrel{d_1^1} \longmapsto\delta_{2p_1}Y_{g/2},\qquad 
Y_g\stackrel{d_1^1} \longmapsto\delta_{2p_1}Y_{g/2}. \]
and
\[ \begin{split}
& X_g\stackrel{d_2^1} \longmapsto\sum_{i=1}^{-p}X_{g/p_i}+\delta_{2p_1}Y_{g/2},\\
& Y_g\stackrel{d_2^1} \longmapsto\sum_{i=1}^{-p}Y_{g/p_i}+\delta_{2p_1}Y_{g/2}.
\end{split} \]
Set $X_m^p=H^2(J, L_{m}^p)$. Let  $X_m^{\bullet}$ be the cochain complex formed by 
$X_m^p$ and $d^1$.  By definition,  for any even positive $q$, $'E_2^{p,q}$ is just the 
$p$-th cohomology group $(X_m^{\bullet},d^1)$.  We calculate the cohomology groups 
one by one:

(1).  $m$ is odd and $d^1=d_1^1$. This is trivial:
\[  ('E_2^{p,q},d_1)={'E}_1^{p,q}=(\Z/2\Z)^{\binom{r}{-p}}.  \]

(2).  $m$ is odd and $d^1=d_2^1$. In this case, if $m=p^n$, it is easy to see that
$H^0(X_{p^n}^{\bullet}, d_2^1)=H^{-1}(X_{p^n}^{\bullet}, d_2^1)=0$. Now if  
$m=m_1m_2$ and $(m_1,m_2)=1$, we can  check 
\[ (X_m^{\bullet}, d_{2m}^1)=(X_{m_1}^{\bullet}, d_{2m_1}^1)\bigotimes (X_{m_2}^{\bullet}, 
d_{2m_2}^1).\]
By K\"unneth's formula,  $H^p(X_m^{\bullet}, d_2^1)=0$. Therefore we have
\[ ('E_2^{p,q},d_2)=\dots=('E_{\infty}^{p,q},d_2)=0.\]

(3).  $m$ is even and $d^1=d_1^1$.  Since $X_m^{\bullet}$ is a $\Z/2\Z$-vector space,
by the formula above about $d_1^1$, we always have
\[ \dim_{\Z/2\Z} \im(X_m^{p}\rightarrow X_m^{p+1})=\binom{r-1}{-p-1}, \]
therefore 
\[ \dim_{\Z/2\Z} \ker(X_m^{p}\rightarrow X_m^{p+1})=2\binom{r}{-p}-\binom{r-1}{-p-1}. \]
Hence
\[ \dim_{\Z/2\Z} H^p(X_{m}^{\bullet}, d_1^1)=2\binom{r}{-p}-\binom{r-1}{-p-1}
-\binom{r-1}{-p}=\binom{r}{-p}. \]
Or  we have
\[  ('E_2^{p,q},d_1)=(\Z/2\Z)^{\binom{r}{-p}}.  \]

(4). $m$ is even and $d^1=d_2^1$.  In this case, if $m=2^k$, 
\[ H^0(X_{2^k}^{\bullet}, d_2^1)=H^{-1}(X_{2^k}^{\bullet}, d_2^1)=\Z/2\Z. \]
Now if $m=2^k m'$, $m'>1$ odd,   set
\[ {X'}_{m'}^p=\bigoplus_{g\mid m'}(<X_g>\bigoplus <Y_g>) \]
and 
\[ d'_2:  X_g\longmapsto\sum_{i=1}^{-p}X_{g/p_i}, \ Y_g\longmapsto
\sum_{i=1}^{-p}Y_{g/p_i}.\]
Then we have 
\[ (X_m^{\bullet}, d^1_{2})=(X_{2^k}^{\bullet}, d^1_{2})\bigotimes (X_{m'}^{'\bullet}, d'_{2}). \]
Similar to the case (2),  we can see $H^p(X_{m'}^{' \bullet}, d'_2)=0$. By K\"unneth's 
formula again,  $('E_2^{p,q},d_2)=H^p(X_{m}^{\bullet}, d_2^1)=0$. 

Combining all the cases above, for $d=d_1$, $q>0$, we have
\begin{equation} \label{e2pq1} ('E_2^{p,q},d_1)=\begin{cases} 
(\Z/2\Z)^{\binom{r}{-p}}, &\text{if}\ q\  \text{odd};\\
0, \ &\text{otherwise} . \end{cases} \end{equation}
For $d=d_2$, $q>0$, we have
\begin{equation} \label{e2pq2} ('E_2^{p,q},d_2)=\begin{cases} 
\Z/2\Z, &\text{if}\ q\ \text{odd}, \ m=2^k, \ p=0\ \text{or} \ -1;\\
0,\ &\text{otherwise} . \end{cases} \end{equation}

Similarly for $F^{\ast, \ast}$,  for $d=d_1$, $q\in \Z$, we have
\[ ('E_{2,F}^{p,q},d_1)=\begin{cases} 
(\Z/2\Z)^{\binom{r}{-p}}, &\text{if}\  q\  \text{odd};\\
0, \ &\text{otherwise} . \end{cases} \]
For $d=d_2$, $q\in \Z$ , we have 
\[ ('E_{2,F}^{p,q},d_2)=\begin{cases} 
\Z/2\Z, &\text{if}\ q\ \text{odd}, \ m=2^k, \ p=0\ \text{or} \ -1;\\
0,\ &\text{otherwise} . \end{cases} \]

 Our next task is to show that $F^{\ast, \ast}$ collapses at $'E_{2,F}^{p,q}$. For this purpose, 
we define a subcomplex $SF^{\ast, \ast}$ of $F^{\ast, \ast}$:
\[ SF^{p, q}= \begin{cases}  (L_m^p,q),\ & \text{if $q$ even} ;\\
(\beta( L_m^p),q)  ,\ & \text{if $q$ odd} .
\end{cases} \]
where
\[ \beta([a,g])= \begin{cases} [a,g],\ & \text{if}\ 2a\neq 0 ;\\
2 [a,g],\ & \text{if}\ 2a= 0 . \end{cases} \]
It is easy to verify $d$ and $\delta$ are well defined in  $SF^{\ast, \ast}$. Moreover,
$H_{\delta}^{q}(SF^{p,\ast})=0$ for every pair $(p,q)$. Therefore the total cohomology
group of $SF^{\ast} $ is trivial.  Hence  to study $F^{\ast,\ast}$, it suffices to study the
double complex $QF^{\ast,\ast}=F^{\ast,\ast}/SF^{\ast,\ast}$. Note that $QF^{\ast,\ast}$
vanishes at the odd rows and at the even rows, $QF^{p,q}$ is nothing but 
$'E_{1,F}^{p,q}$ we just got above.  Moreover, the induced differential $d$ in $QF^{\ast,\ast}$ is
 nothing but $d^1$ in $'E_{1,F}^{\ast,\ast}$.  On one hand
\[ H^n(QF^{\ast}) =H^n(F^{\ast}) =\hat{H}^n(J, H_d^0(L_m))=\bigoplus_{q} \, 
'E_{\infty,F}^{n-q,q}; \]
on the other hand, we have
\[ \begin{split} H^n(QF^{\ast})=&\frac
{\ker\, (d+\delta: \bigoplus_{q} QF^{n-q, q}\rightarrow  \bigoplus_{q} QF^{n+1-q, q})}
{\im\, (d+\delta: \bigoplus_{n} QF^{n-1-q, q}\rightarrow \bigoplus_{q} QF^{n-q, q})} \\ 
=& \bigoplus_{q} \frac{\ker\, (d: QF^{n-q, q}\rightarrow QF^{n+1-q, q})}
{\im\, (d: QF^{n-1-q, q}\rightarrow QF^{n-q, q})}\\
=& \bigoplus_{q} \, 'E_{2,F}^{n-q,q}. \end{split} \]
Therefore for any pair $(p,q)$, 
\[ 'E_{2,F}^{p,q}='E_{\infty,F}^{p,q}. \]
In particular, for $n=1$,
\begin{equation} \label{cru}
 \hat{H}^1(J, H_d^0(L_m))=\bigoplus_{q}\, 'E_{2,F}^{1-q,q}=\bigoplus_{q}\,  'E_{2}^{1-q,q}. 
\end{equation}
By results of  \eqref{e2pq1} and \eqref{e2pq2}, we easily have
\begin{theo} \label{coho} 
The group $J=\{1,c\}$ acts trivially on the cohomology groups $H^i(J, \co_m)$
and $H^i(J, U_m)$ for $i=1$ or $2$, moreover, 
\begin{equation} H^1(J, \co_m)=H^2(J, \co_m)=\begin{cases} \Z/2\Z, \ & \text{if}\  
m=2^k;\\0, \ & \text{otherwise}.
\end{cases} \end{equation}
and 
\begin{equation} H^1(J, U_m)=H^2(J, U_m)=(\Z/2\Z)^{2^{r-1}}.  \end{equation}
\end{theo}
\begin{rem} 1. The second statement is first proved in Yamamoto~\cite{Yamamoto}.
The spectral sequence method employed here makes the calculation significantly simpler 
than those in \cite{Yamamoto} and in \cite{S1}. Moreover, this same spectral sequence
method can also be applied  to the universal distribution of 
higher rank, thus recover the results in Kubert~\cite{Kubert2}. 

2. For any cyclic group $C\in G$ which has trivial intersection with $G_{p^i}$ for $p^i\| m$, 
we can also obtain similar result without any extra difficult. 
 \end{rem}
\begin{prop} \label{index}
\[ I(L_m, d_1;\theta)=\begin{cases} 2,\ &  if\ r=1;\\ 2^{2^{r-2}}, \ &\text{if}\ r>1. \end{cases}
\qquad  I(L_m, d_2;\theta)=\begin{cases} 2,\ &  if\ m=2^k;\\1,\ & otherwise. \end{cases} \] 
\end{prop}
\begin{proof} By the identity~\eqref{cru},  the condition in Proposition~\ref{lapse}
is satisfied.  For $d=d_1$, by \eqref{e2pq1}, then  the exponent of $2$ in $I(L_m, d_1;\theta)$ is
equal to 
\[ \begin{split} \sum_{\substack{p+q\leq 0\\q>0\ odd}}(-1)^{p+1} \binom{r}{-p}=
&\sum_{p=-r}^{-1}\sum_{\substack{q\leq -p\\q\ odd}}(-1)^{p+1} \binom{r}{-p}\\
=&\sum_{p=1}^{r}\left [\frac{p+1}{2}\right ] (-1)^p \binom{r}{p}\\
=&\sum_{k} k\binom{r}{2k-1}-\sum_{k} k\binom{r}{2k}\\
=& \begin{cases} 1,\ &if \ r=1;\\2^{r-2}, \ &if \ r>1. \end{cases} \end{split} \] 
The case $d=d_2$ immediately follows from Proposition~\ref{lapse} and   \eqref{e2pq2}.
\end{proof}

\section{Sinnott's index formula}
In this section we set to prove the following theorem:
\begin{theo}[See \cite{S1}, Theorem] Let $R=\Z[G_m]$ and let $S$ be the Stickelberger 
ideal of $\Q(\zeta_m)$. Then
\[ [R^{-}:S^{-}]=2^a h^{-}, \]
where $a=0$ if $r=1$ and $a=2^{r-2}-1$ if $r>1$.
\end{theo}
\begin{note}In this section,  the subscript $m$ is omitted from our notations(i.e., 
$G$ is the Galois group $G_m$ and so on). $p$ is always regarded as a prime factor of $m$. The 
superscript ``-'', is in accordance with the superscript ``$\theta$'' in the previous sections. 
\end{note}
\begin{proof} We consider the following diagram:
\[  \xymatrix{
{\R U^{-}}\ar[r]^{\varphi^{-}}\ar[dr]_{\alpha^{(s)}} & 
{\R \co^{-}}\ar[r]^{\psi^{(s)}}&{\R \co^{-}}\ar[dl]^{\beta}\\
&{\R[G]^{-}}&
 } \]
where $s>1$ and $\varphi^{-}=\varphi|_{\R U^{-}}$,
\[ \psi^{(s)}([x]-[-x])=\sum_{(n,m)=1}\frac{[nx]-[-nx]}{n^s}, \]
\[ \beta([x]-[-x])=\frac{1}{2\pi i}\sum_{t\in (\Z/m\Z)^{\times}} (\exp({2\pi i xt})-
\exp(-{2\pi i xt})) \sigma^{-1}_t \]
and
\[ \alpha^{(s)}=\beta\circ \psi^{(s)}\circ \varphi^{-}. \]
$\psi^{(s)}$ is well defined and all the above maps are isomorphisms of vector spaces. Then 
we have
\begin{equation} \label{sin1} \begin{split}
&(R^{-}: \alpha^{(s)}(U^{-}))\\
=&(R^{-}:\beta(\co^{-})) \cdot
(\beta(\co^{-}): \beta\psi^{(s)}(\co^{-}))\cdot
(\beta\psi^{(s)}(\co^{-}): \alpha^{(s)}(U^{-}))\\
=&(R^{-}:\beta(\co^{-}))\cdot
(\co^{-}: \psi^{(s)}(\co^{-}))\cdot
(\co^{-}: \varphi^{-}(U^{-}))
\end{split} \end{equation}
Here for the second equality, we use the property that if $V_1$ and $V_2$ are two 
vector spaces and $f$ is an isomorphism from $V_1$ to $V_2$, then $(A:B)_{V_1}
=(f(A):f(B))_{V_2}$. Now for the three factors at the last line of \eqref{sin1}, we have:

\begin{la} \label{pqrf} \begin{equation} \label{sin2}
(R^{-}:\beta(\co^{-}))=\begin{cases} (2\pi)^{-\varphi(m)/2}\sqrt{d(K_m)/d(K^{+}_m)},
\ &\text{if}\ m\neq 2^k;\\ \frac{1}{2}(2\pi)^{-\varphi(m)/2}\sqrt{d(K_m)/d(K^{+}_m)},
\ &\text{if}\ m=2^k . \end{cases}
\end{equation}  \end{la}
\begin{proof}[Proof of Lemma~\ref{pqrf}]  We first consider the following diagram
which is exact at the rows:
\[
\begin{CD}
0@>>> \co^{+} @>{i}>> \co @>{1-c}>>\im (1-c) @>>>0\\
@. @AA{id}A @AA{i}A  @AA{i}A @.\\
0@>>> \co^{+} @>{i}>> \co^{+} \oplus\co^{-}@>{1-c}>>2\co^{-} @>>>0
\end{CD}  \]
\vskip 0.4cm
\noindent where $i$ is the natural inclusion map.  By Theorem~\ref{coho}, if 
$m$ is not a power of $2$, $\co^{-}=\im (1-c)$; if $m$ is a power of $2$, then 
$\co^{-}/\im (1-c)=\Z/2\Z$. Therefore,
\[ (\co:\co^{+} \oplus\co^{-})
=(\im(1-c):2\co^{-})=\begin{cases} 2^{\varphi(m)/2},\ &\text{if}\ m\neq 2^k;\\
2^{\varphi(m)/2-1},\ &\text{if}\ m= 2^k. \end{cases}  \]
Now let $T$ be the map from $\C\co$ to $\C [G]$ such that $T([x])=\sum_{t} 
\exp(2\pi i tx) \sigma^{-}_t$, then we have $T|_{\C\co^{-}}=2\pi i\beta|_{\C \co^{-}}$.
Then on one hand,  
\[ (R^{+}\oplus R^{-}:T(\co^{+}\oplus\co^{-}))=(R^{+}:T(\co^{+}))\cdot (R^{-}:T(\co^{-})) \]
on the other hand, 
\[ (R^{+}\oplus R^{-}:T(\co^{+}\oplus\co^{-}))=(R^{+}\oplus R^{-}:R)\cdot
(R:T(\co))\cdot (\co: \co^{+}\oplus\co^{-}). \]
But we know $(R^{+}\oplus R^{-}:R)=2^{-\varphi(m)/2}$,  and by the definition of $T$, $(R:T(\co))=
\sqrt{d(K)}$ and $(R^{+}:T(\co^{+}))=\sqrt{d(K^{+})}$. Now the lemma follows from the above results 
and
\[  (R^{-}:\beta(\co^{-}))=(2\pi)^{-\varphi(m)/2} (R^{-}: 2\pi i\beta(\co^{-})). \]
\renewcommand{\qed}{} \end{proof}

\begin{la} \label{mnn} Let $S=\{p:\, p\mid m\}$, then
\begin{equation} \label{sin3}
(\co^{-}: \psi^{(s)}(\co^{-}))
=\prod_{\chi\ odd} L_S(s,\chi). 
\end{equation}\end{la}
\begin{proof}[Proof of Lemma~\ref{mnn}] Note that if we let
\[ \Theta_{S}(s)=\sum_{(n.m)=1}\frac{\sigma_n}{n^s}, \]
then $\psi^{(s)}$ is just the left multiplication of $\Theta_{S}(s)$ on $\R \co^{-}$. 
By \cite{S2} Lemma 1.2(b), we have
\[ (\co^{-}: \psi^{(s)}(\co^{-}))=\prod_{\chi\ odd}
\chi( \Theta_{S}(s))=\prod_{\chi\ odd} L_S(s,\chi). \] 
\renewcommand{\qed}{} \end{proof}

\begin{la} \label{mnl}
\begin{equation} \label{sin4} 
(\co^{-}: \varphi^{-} (U^{-}))
=\begin{cases} 2^{-2^{r-2}}
\underset{p\mid m}{\prod}\ \underset{\chi\ odd}{\prod} (1-\chi(p)^{-1})p^{-1},\ & if\ r>1;\\
\frac{1}{2},\ &if\ r=1, \ p\neq 2;\\
1,\ &if\ m=2^k. \end{cases} 
\end{equation}\end{la}
\begin{proof}[Proof of Lemma~\ref{mnl}] This follows from the abstract index formula~\eqref{aif} ,
Proposition~\ref{deter1}  and Proposition~\ref{index}. 
\renewcommand{\qed}{} \end{proof}

Now let $s$ approach $1$, then 
\begin{equation} \label{sin5} \begin{split}
  \lim_{s\rightarrow 1} \alpha^{(s)}([x]-[-x])= &
\lim_{s\rightarrow 1} \beta\psi^{(s)}H^0(\varphi)([x]-[-x])\\=& \frac{1}{2\pi i}
\sum_{t} \sigma^{-1}_t\sum_{n\in \N}\frac{\exp(2n\pi i xt)-\exp(-2n\pi i xt)}{n}\\
=& \sum_t (\frac{1}{2}-\{xt\})\sigma^{-1}_t.
\end{split} \end{equation}
If we let $\alpha= \lim_{s\rightarrow 1} \alpha^{(s)}$, by \eqref{sin1}, \eqref{sin2},\eqref{sin3} 
and \eqref{sin4}, with the class number formula,
\[ h^{-} =(2\pi)^{-\varphi(m)/2} \prod_{\chi\ odd} L(1,\chi) \sqrt{d(K_m)/d(K^{+}_m)}\,\omega\, Q, \]
and since $(U^{-}: (1-c)U)=2^{2^{r-1}}$, then we have 
\begin{equation} \label{sin6} \begin{split}
(R^{-}: \alpha((1-c)U))=&\lim_{s\rightarrow 1}
(R^{-}: \alpha^{(s)}(U^{-}))\cdot (U^{-}: (1-c)U)\\
=&\begin{cases} \frac{h^{-}}{\omega Q}\cdot 2^{2^{r-2}}, \ &if \ r>1;\\
\frac{h^{-}}{\omega Q}, \ &if \ r=1. \end{cases} \end{split}
\end{equation}
But  by \eqref{sin5}, 
$\alpha((1-c)U)$ is nothing but $e^{-}S'$ in \cite{S1}. and by \cite{S1}, Lemma 3.1,
we have $(e^{-}S':S^{-})=\omega$. This is enough to finish the proof of the theorem..
\end{proof}

\end{document}